\documentclass{amsart}

\usepackage{amsfonts,amssymb,amscd,amsmath,enumerate,url,verbatim}
\usepackage[all]{xy}
\SelectTips{eu}{}
\xyoption{curve}

\newcommand{\lra}{\longrightarrow}

\newcommand{\xra}{\xrightarrow}

\newcommand{\cent}[1]{{{#1}^{\mathsf c}\!}}

\newcommand{\ov}{\bar}

\newcommand{\wh}{\widehat}
\newcommand{\wt}{\widetilde}

\newcommand{\ges}{{\scriptscriptstyle\geqslant}}
\newcommand{\sls}{{\scriptscriptstyle<}}

\newcommand{\dcat}[1]{{\mathsf{D}(#1)}}
\newcommand{\lotimes}[1]{\otimes^{\bf L}_{#1}}

\newcommand{\shift}{{\sf\Sigma}}
\newcommand\col{\colon}
\newcommand{\cone}[1]{{\mathsf{cone}}(#1)}
\newcommand\dd{\partial}
\newcommand{\hh}[1]{\operatorname{H}(#1)}
\newcommand{\hch}[3]{\operatorname{H}^{#1}(#3\var #2)}
\newcommand{\HH}[2]{\operatorname{H}_{#1}(#2)}
\newcommand{\CH}[2]{\operatorname{H}^{#1}(#2)}

\newcommand{\Ext}[4]{\operatorname{Ext}^{#1}_{#2}(#3,#4)}
\newcommand{\Hom}[3]{\operatorname{Hom}_{#1}(#2,#3)}

\newcommand\Ker{\operatorname{Ker}}

\newcommand\spec{\operatorname{Spec}}

\newcommand{\ann}{\operatorname{ann}}

\newcommand\depth{\operatorname{depth}}

\newcommand{\rank}{\operatorname{rank}}

\newcommand\idmap{\operatorname{id}}

\newcommand\pd{\operatorname{proj\,dim}}

\newcommand\fm{{\mathfrak m}}

\newcommand\fp{{\mathfrak p}}
\newcommand\fq{{\mathfrak q}}

\newcommand\bs{\boldsymbol}

\newcommand\bsf{{\boldsymbol f}}

\newcommand\bsx{{\boldsymbol x}}

\newcommand\bsz{{\boldsymbol z}}

\newcommand\BZ{{\mathbb Z}}

\newcommand\sse{{\mathsf e}}

\newcommand{\eps}{{\varepsilon}}
\newcommand{\vf}{{\varphi}}

\newcommand{\var}{{\hskip1pt\vert\hskip1pt}}

\theoremstyle{plain}
\newtheorem{theorem}{Theorem}[section]

\newtheorem{etheorem}[theorem]{Existence Theorem}
\newtheorem{proposition}[theorem]{Proposition}
\newtheorem{lemma}[theorem]{Lemma}

\newtheorem{itheorem}{Theorem}{\Alph{theorem}}

\theoremstyle{definition}

\newtheorem{chunk}[theorem]{}

\theoremstyle{remark}
\newtheorem{example}[theorem]{Example}

\newenvironment{bfchunk}{\begin{chunk}\textit}{\end{chunk}}

\numberwithin{equation}{theorem}

\newcommand{\Supp}[1][\mca]{\operatorname{Supp}_{#1}}
\newcommand{\msup}[1][\mca]{\operatorname{Supp}^{\scriptscriptstyle +} 
_{#1}}
\newcommand{\hsup}[2][\mca]{\operatorname{Supp}^*_{#1}(#2)}
\newcommand{\hhsup}[3][\mca]{\operatorname{Supp}^*_{#1}(#2,#3)}
\newcommand{\mvar}[2][\mca]{V_{\!#1}(#2)}
\newcommand{\hhvar}[2][\mca]{\displaystyle V^{*}_{\!#1}(#2)}

\newcommand{\Max}[1]{\operatorname{Max}{#1}}

\newcommand{\thick}[2]{{\mathsf{Thick}}_{#1}(#2)}

\newcommand{\proj}{\operatorname{Proj}}
\newcommand{\syz}[3][R]{{\Omega^{#1}_{#2}(#3)}}
\newcommand{\koz}[2]{#2/\!\!/#1}

\newcommand{\mca}{{\mathcal{A}}}

\newcommand{\mci}{{\mathcal{I}}}
\newcommand{\mcl}{{\mathcal{L}}}
\newcommand{\mcm}{{\mathcal{M}}}
\newcommand{\mcn}{{\mathcal{N}}}
\newcommand{\mcr}{{\mathcal{R}}}

\begin{document}
\title[Cohomological supports]{Constructing modules with prescribed \ 
\ cohomological support}
\date{30th July 2007}

\author[L.~L.~Avramov]{Luchezar L.~Avramov}
\address{Department of Mathematics, University of Nebraska, Lincoln,  
NE 68588, U.S.A.}
\email{avramov@math.unl.edu}

\author[S.~B.~Iyengar]{Srikanth B.~Iyengar} \address{Department of
Mathematics, University of Nebraska, Lincoln, NE 68588, U.S.A.}
\email{iyengar@math.unl.edu}
\thanks{Research partly supported by NSF grants DMS 0201904 (L.L.A),  
DMS 0602498 (S.I.)}

   \begin{abstract}
A cohomological support, $\operatorname{Supp}^*_{\mathcal A}(M)$,
is defined for finitely generated modules $M$ over a left noetherian
ring $R$, with respect to a ring $\mathcal A$ of central cohomology
operations on the derived category of $R$-modules. It is proved that
if the $\mathcal A$-module $\operatorname{Ext}^*_R(M,M)$ is noetherian
and $\operatorname{Ext}^*_R(M,R)=0$ for $i\gg0$, then every closed
subset of $\operatorname{Supp}^*_{\mathcal A}(M)$ is the support of
some finitely generated $R$-module.  This theorem specializes to known
realizability results for varieties of modules over group algebras, over
local complete intersections, and over finite dimensional algebras over
a field.  The theorem is also used to produce large families of finitely
generated modules of finite projective dimension over commutative local
noetherian rings.
   \end{abstract}

\dedicatory{To Phil Griffith, algebraist and friend.}



\maketitle

\tableofcontents

\section*{Introduction}

Quillen introduced methods from algebraic geometry to the study of
cohomology rings of finite groups in a seminal paper, \cite{Qu}.
His ideas and techniques have led to the appearance of a number of
highly developed theories, which provide insight into the structure of
an algebraic object through some geometric `variety' attached to it.
Use of such geometric invariants has been crucial to progress on a  
number
of difficult problems.

Variety theories share certain formal properties needed in applications.
Some of them guarantee that homologically similar modules, such as all
syzygy modules of a given module, have the same variety.  Modules with
distinct varieties are therefore expected to exhibit quantifiable
differences in homological behavior.  For this reason, a description
of all the varieties produced by a given theory is a useful tool for
classifying homological patterns.

The prototype theory applies to all finite dimensional representations
of a finite group; see \cite{Be:rc2} for a detailed exposition.  It has
been extended to representations of finite dimensional cocommutative
Hopf algebras, \cite{FS,SFB}.  Parallel theories have been constructed
for finitely generated modules over finite dimensional self-injective
algebras, \cite{EHSST, SS}, and over local complete intersection rings,
\cite{Av:vpd, AB}.  Historically, in each concrete case the proofs of
the formal properties of a theory and of the relevant realizability
theorem have involved delicate arguments specific to that context.

We are interested in modules over a fixed associative ring $R$.

The vehicle for passing from algebra to geometry is provided by a choice
of commutative graded ring $\mca$ of central cohomology operations
on the derived category of $R$.  In the examples above there are
natural candidates for $\mca$: the even cohomology ring of a group
(or a Hopf algebra); the even subalgebra of the Hochschild cohomology
of an associative algebras; the polynomial ring of Gulliksen operators
over a complete intersection.  However, other choices are possible and
sometimes are desirable.

For each pair $(M,N)$ of $R$-modules the graded group $\Ext*RMN$ has
a natural structure of graded $\mca$-module.  The set
   \[
\hhsup MN = \{\fp\in \proj \mca \mid \Ext *RMN_\fp\ne 0\}
   \]
where $\proj\mca$ is the space of all essential homogeneous prime ideals
in $\mca$ with the Zariski topology, is called the \emph{cohomological
support} of $(M,N)$. The cohomological support of $M$ is the set
$\hhsup MM$.

The principal contribution of this work is a method for constructing
modules with prescribed cohomological supports.  Part of our main
result reads:

   \begin{itheorem}
  \label{itheorem1}
Let $R$ be a noetherian ring and let $M$ and $N$ be finite $R$-modules,
such that the graded $\mca$-module $\Ext *RMN$ is noetherian.

If $\Ext iRMR=0$ holds for all $i\gg0$, then for every closed subset $X$
of $\hhsup MN$ there exist finite $R$-modules  $M_X$ and $N_X$ such that
   \[
\hhsup{M_X}N=X=\hhsup M{N_X}\,.
   \]
Moreover, when $N=M$ one can choose $N_X=M_X$.
  \end{itheorem}

Suitable specializations of Theorem \ref{itheorem1} yield several known
realizibility results: See Section \ref{bialgebras} for Hopf algebras
and Section \ref{Associative algebras} for associative algebras.
Their earlier proofs were modeled on Carlson's Tensor Product Theorem
\cite{Ca} for varieties over group algebras; they rely heavily on the
nature of $R$ (in the first case) or on that of $\mca$ (in the second).

Theorem \ref{itheorem1}  is proved in Section \ref{realizability
by modules}, based on work in Sections \ref{syzygy complexes} and
\ref{cohomological supports}.  Our argument requires few structural
restrictions on $R$ and none on $\mca$ itself.  The crucial input is
the noetherian property of $\Ext *RMN$ as a module over $\mca$.

Another application of Theorem \ref{itheorem1} goes into a completely
different direction:

   \begin{itheorem}
   \label{itheorem2}
Let $(Q,\fq,k)$ be a commutative noetherian local ring and ${\bsf}$
a $Q$-regular sequence of length $c$ contained in $\fq^2$.

For $R=Q/Q\bsf$ and $\ov k$ an algebraic closure of $k$ there exists
a map
   \[
V\col
\left\{\begin{gathered}
\text{isomorphism classes $[M]$ }
\\
\text{of finite $R$-modules with }
\\
\text{$\pd_QM<\infty$}
\end{gathered}\right\}
\lra
\left\{\begin{gathered}
\text{closed algebraic }
\\
\text{sets}\ \ X\subseteq{\mathbb P}^{c-1}_{\ov k}
\\
\text{defined over }k
\end{gathered}\right\}
   \]
with the following properties:
   \begin{enumerate}[\quad\rm(1)]
  \item
$V$ is surjective.
  \item
$V([M])=\varnothing$ if and only if $\pd_RM<\infty$.
  \item
$V([M])=V([\Omega^R_n(M)])$ for every syzygy module $\Omega^R_n(M)$.
  \item
$V([M])=V([M/\bsx M])$ for every $M$-regular sequence $\bsx$ in $R$.
   \end{enumerate}
  \end{itheorem}

This result is surprising.  Indeed, it exhibits large families of  
modules
of finite projective dimension over any ring $Q$ with $\depth Q\ge2$,
contrary to a commonly held perception that finite projective dimension
is `rare' over singular commutative rings.  Furthermore, the remaining
statements ascertain that modules mapping to distinct closed cones in
${\ov k}^c$ cannot be linked by any sequence of standard operations
known to preserve finite projective dimension.

In Section \ref{local complete intersections} we prove Theorem
\ref{itheorem2}, and deduce from it a recent theorem on the existence
of cohomological varieties for modules over complete intersection local
rings.  For the latter we establish a descent result of independent
interest.

In this paper varieties of modules are discussed in the broader context
of varieties of complexes.  The resulting marginal technical  
complications
are easily offset by a gain in flexibility:  We first realize a given  
set
as the cohomological support of a bounded complex by using constructions
whose effect is easy to track.  To show that this set is also the  
support
of a module we use `syzygy complexes', a notion introduced and discussed
in Section \ref{syzygy complexes}.

This paper is part of an ongoing study of cohomological supports
of modules over general associative rings.  In \cite{AI} we focus on
proving existence of variety theories with desirable properties under a
small set of conditions on a ring, its module(s), and a ring of central
cohomological operators.  The properties that have to be established are
clarified in \cite{BIK} by Benson, Iyengar, and Krause, who investigate
a notion of support for triangulated categories equipped with an action
by a central ring of operators. On the other hand, the methods of this
paper can be adapted to prove realizability results in that context.
Of particular interest is the case of certain monoidal categories, where
work of Suarez-Alvarez, \cite{Su} provides natural candidates for rings
of operators.

\section{Syzygy complexes}
\label{syzygy complexes}

In this section we recall a few basic concepts of DG homological  
algebra,
following \cite{AFH}, and extend the notion of syzygy from modules
to complexes.

Let $R$ be an associative ring and $\dcat R$ the full derived  
category of
left $R$-modules.  We write $\simeq$ to indicate a quasi-isomorphism of
complexes; these are the isomorphisms in $\dcat R$.  The symbol $\cong$
is reserved for isomorphisms of complexes, and $\equiv$ is used to  
denote
homotopy equivalences. Given a complex $M$ of $\dcat R$, we write $ 
\thick
RM$ for its \emph{thick closure}, that is to say, the intersection of
the thick subcategories of $\dcat R$ containing $M$.

\begin{bfchunk}{Semiprojective complexes.}
   \label{chunk:semiproj}
A complex $P$ of $R$-modules is called \emph{semiprojective} if $\Hom
RP-$ preserves surjective quasi-isomorphisms; equivalently, if $P$ is a
complex of projective modules and $\Hom RP-$ preserves quasi- 
isomorphisms.
The following properties are used in the proofs below.

Every quasi-isomorphism of semiprojective complexes is a homotopy
equivalence.  Every surjective quasi-isomorphism to a semiprojective
complex has a left inverse.  Every semiprojective complex $C$ with $\hh
C=0$ is equal to $\cone{\idmap^B}$ for some complex $B$ of projective
modules with zero differential.
   \end{bfchunk}

   \begin{lemma} \label{lem:semiproj-complexes}
If $\pi\col P\to Q$ is a quasi-isomorphism of semiprojective complexes
of $R$-modules and $n$ is an integer, then there is a homotopy  
equivalence
   \[
P_{\ges n}\oplus\shift^n Q'\equiv Q_{\ges n}\oplus \shift^n P'
   \]
where $P'$ and $Q'$ are projective $R$-modules.
     \end{lemma}

\begin{proof}
Assume first that $\pi$ is surjective.  It then has a left inverse,  
hence
one gets $P\cong Q\oplus E$ with $E=\Ker(\pi)$.  This implies that $E$
is semiprojective with $\hh E=0$, and hence $E=\cone{\idmap^{F}}$ for
some complex $F$ of projective $R$-modules with $\dd^{F}=0$.  Hence one
gets a quasi-isomorphism
   \[
P_{\ges n}\cong Q_{\ges n}\oplus\cone{\idmap^{F_{\ges n}}}\oplus\shift^n
F_{n-1}\,.
   \]
The canonical map $P_{\ges n}\to Q_{\ges n}\oplus\shift^n F_{n-1}$
is thus a homotopy equivalence, as  $\cone{\idmap^{F_{\ges n}}}$ is
homotopy equivalent to $0$. This settles the surjective case.

In general, $\pi$ factors as $P \to \wt P \xra{\,\psi\,}Q$, where $ 
\wt P$
is equal to $P\oplus\shift^{-1}\cone{\idmap^Q}$ and $\psi$ is the sum of
$\pi$ and the canonical surjection $\shift^{-1}\cone{\idmap^Q}\to Q$.
Thus, $\psi$ is a surjective quasi-isomorphism of semi-projective
complexes.  So is the canonical map $\wt P\to P$.  The already settled
case yields homotopy equivalences
   \[
P_{\ges n}\oplus\shift^n Q'\gets \wt P\to Q_{\ges n}\oplus \shift^n  
P'\,.
   \]
for appropriate projective modules $P'$ and $Q'$.
  \end{proof}

\begin{bfchunk}{Syzygy complexes.}
   \label{chunk:syzygies}
Let $M$ be a complex of $R$-modules.

A \emph{semiprojective resolution} of $M$ is a quasi-isomorphism $P 
\to M$
from a semiprojective complex $P$.  Every complex $M$ has one, and it
is unique up to homotopy equivalence.  Thus, the preceding result may
be viewed as a homotopical version of Shanuel's Lemma.  Based on it,
we introduce a homotopical version of the notion of syzygy module.

For each $n\in\BZ$ let $\syz nM$ stand for any complex
$\shift^{-n}(P_{\ges n})$, where $P$ is a semiprojective resolution
of $M$, and call it an $n$th \textit{syzygy complex} of $M$ over $R$.
Its dependence on the choice of $P$ is made precise by the preceding
lemma.

Any complex $P$ of projective modules with $P_i=0$ for $i\ll0$ is
semiprojective.  Thus, when $M$ is an $R$-module and $P$ is its  
projective
resolution the complex $\shift^{-n}(P_{\ges n})$ is isomorphic in $\dcat
R$ to an $n$th syzygy module of $M$.
   \end{bfchunk}

The next  lemma expands upon the last observation.

  \begin{lemma}
   \label{lem:syzygies}
If $M$ is a complex of $R$-modules,  $s=\sup\{i\var \HH iM\ne 0\}$,
and $n$ is an iteger with $n\ge s$, then $\syz nM$ is quasi-isomorphic
to $\HH0{\syz nM}$.
     \end{lemma}

\begin{proof}
Let $P\to N$ be a semiprojective resolution with $\syz
nM=\shift^{-n}(P_{\ges n})$.  For $i\ge n+1$ one has isomorphisms $\HH
i{P_{\ges n}}\cong \HH iP\cong \HH iM=0$, the first one of which comes
from the exact sequence of complexes
   \begin{equation}
     \label{eq:trunk}
0\to P_{\sls n}\to P\to P_{\ges n}\to 0\,.  \qedhere
   \end{equation}
  \end{proof}

  \begin{bfchunk}{Cohomology.}
    \label{chunk:cohomology}
Let $M$ be a complex of $R$-modules and $P\to M$ a semiprojective
resolution.  For every complex $N$ and each $i\in\BZ$ the abelian group
\[
\Ext iRMN= \HH {-i}{\Hom RPN}=\CH{i}{\Hom RPN}
   \]
is independent of the choice of resolution $P$, see \ref 
{chunk:semiproj};
it is a module over $R^{\mathsf c}$, the center of the ring $R$.
For modules $M$, $N$ this is the usual gadget, see \ref{chunk:syzygies}.
   \end{bfchunk}

Over noetherian rings syzygy modules inherit finiteness properties of
the original module.  We show that syzygy complexes behave similarly.

\begin{lemma}
\label{lem:thick}
If $R$ is a noetherian ring and $M$ is a complex with $\hh M$ a finite
$R$-module, then one can find a syzygy complex $\syz nM$ in
$\thick R{M\oplus R}$.

Furthermore, for every complex $C\in\thick R{M\oplus R}$ the  
following hold.
\begin{enumerate}[\quad\rm(1)]
\item The $R$-module $\hh C$ is noetherian.
\item $\Ext{\gg0}RMN=0$ for a bounded complex $N$ implies
$\Ext{\gg0}RCN=0$.
\item $\Ext{\gg0}RMR\!=\!0$ implies $\Ext{\gg0}RCF=0$ for every
projective $R$-module $F$.
\end{enumerate}
\end{lemma}

\begin{proof}
Under the hypotheses on $R$ and $M$, one can choose a semiprojective
resolution $P\simeq M$ with each $P_i$ finite and $P_i=0$ for $i\ll 0$.
It follows that $P_{<n}$ is in $\thick RR$, so the exact sequence
\eqref{eq:trunk} yields $P_{\ges n}\in\thick R{M\oplus R}$.

The complexes $L$ with $\hh L$ finite form a thick subcategory of
$\dcat R$.  As it contains $M$ and $R$, it contains  $\thick R{M\oplus
R}$ as well.  This proves (1).  A similar argument settles (2).  For $M$
and $F$ as in (3) there is an isomorphism
\[
\Ext *RMF \cong \Ext *RMR\otimes_RF\,,
\]
which one can get by using the resolution $P$ above.  Thus,
$\Ext{\gg0}RMR=0$ implies $\Ext{\gg0}RMF=0$.  Now (2) yields
$\Ext{\gg0}RCF=0$, as desired.
   \end{proof}

\section{Graded rings}
\label{graded rings}

Here we describe notation and terminology for dealing with graded  
objects.

\begin{bfchunk}{Graded modules.}
\label{chunk:noetherian}
Let $\mca$ be a commutative ring that is \textit{non-negatively graded}:
$\mca=\bigoplus_{i\in\BZ}\mca^i$ with $\mca^i\mca^j\subseteq\mca^{i+j}$
and $\mca^i=0$ for $i<0$.

Modules over $\mca$ are $\BZ$-graded: $\mcm= \bigoplus_{j\in\BZ}\mcm^j$
with $\mca^i\mcm^j\subseteq\mcm^{i+j}$.  For such an $\mcm$ \emph 
{finite}
means finitely generated, \emph{eventually noetherian} means $\mcm^{\ges
j}$ is noetherian for $j\gg0$, and \emph{eventually zero} means
$\mcm^{\ges j}=0$ for $j\gg0$.

The \emph{annihilator} of $\mcm$ is the set $\ann_{\mca}\mcm =
\{a\in \mca\mid a\mcm=0\}$.  It is a homogeneous ideal in $\mca$, so
$\mca/\ann_{\mca}\mcm$ is a graded ring and $\mcm$ is a graded module
over it.  When $\mcm$ is noetherian so is $\mca/\ann_{\mca}\mcm$, so
modulo $\ann_{\mca}\mcm$ every ideal in $\mca$ is generated by finitely
many homogeneous elements.
  \end{bfchunk}

\begin{bfchunk}{Supports.}
   \label{chunk:supp}
Let $\spec\mca$ be the space of prime ideals of $\mca$, with the
Zariski topology.  For an $\mca$-module $\mcm$, set
   \begin{align*}
\Supp{\mcm}&= \{\fp\in\spec{\mca}\mid\mcm_\fp\ne0\}\,;
   \\
\proj{\mca}&=\{\fp\in\spec\mca\mid\fp
   \text{ homogeneous and }
\fp\not\supseteq \mca^{\ges1}\}\,;
   \\
\msup[\mca]{\mcm}&= \Supp{\mcm}\cap\proj{\mca}\,.
   \end{align*}

The following properties of graded $\mca$-modules $\mcl,\mcm$, and $ 
\mcn$
follow from the definition of support and the exactness of localization.
\begin{enumerate}[{\quad\rm(1)}]
  \item
If $\mcl \xra{\iota} \mcm\xra{\eps}\mcn$ is an exact sequence, then
\[
  \msup \mcm \subseteq \msup \mcl\,\cup\,\msup\mcn\,;
\]
equality holds when $\iota$ is injective and $\eps$ is surjective.
  \item
For each $i\in\BZ$, one has
\(
\msup {(\mcm^{\ges i})}=\msup \mcm\,.
\)
  \item
If some $\mcm^{\ges n}$ is finite, say, if $\mcm$ is eventually
noetherian, then
\[
\msup \mcm=\{\fp\in \proj\mca\mid \fp\supseteq \ann_\mca(\mcm^{\ges  
i})\}
\]
holds for every $i\ge n$; thus, $\msup \mcm$ is a closed subset of
$\proj\mca$.
  \item
If the $\mca$-modules $\mcm$ and $\mcn$ are finite, then
  \[
\msup(\mcm\otimes_\mca\mcn)=\msup\mcm\,  \cap \, \msup\mcn\,.
  \]
\item If $\mcm$ is eventually zero, then $\msup\mcm=\varnothing$.
The converse holds when $\mcm$ is eventually noetherian over $\mca$.
\end{enumerate} \end{bfchunk}

In some cases, supports have a natural geometric interpretation.

\begin{bfchunk}{Varieties.}
\label{chunk:var}
Let $k$ be a field and $\ov k$ an algebraic closure of $k$.  Assume that
the graded ring $\mca$ has $\mca^0=k$ and is generated over $k$ by
finitely many homogeneous elements of positive degree.  For each graded
$\mca$-module $\mcm$ set
   \[
\mvar {\mcm} = \left(\Supp[\ov\mca]{(\mcm\otimes_k\ov
k)}\cap\Max{\ov\mca}\right) \cup\{{\ov\mca}^{\ges1}\}
   \]
where $\ov\mca$ denotes the ring $\mca\otimes_k\ov k$ and $\Max{\ov 
\mca}$
the set of its maximal ideals.

Let $\mcm$ be a finite graded $\mca$-module.  The subset $\mvar {\mcm}$
of $\Max{\ov\mca}$ then is closed in the Zariski topology; it is
also $k$-rational and conical, in the sense that it can be defined by
homogeneous elements in $\mca$.   The Nullstellensatz implies that each
one of the sets $\mvar{\mcm}$ and $\msup{\mcm}$ determines the other.
   \end{bfchunk}

The graded rings and modules of interest in this paper are generated by
cohomological constructions, which we recall below.

  \begin{bfchunk}{Products in cohomology.}
    \label{chunk:products in cohomology}
Let $M$ and $N$ be complexes of $R$-modules, and let $P\to M$ and $Q 
\to N$
be semiprojective resolutions.  For each $i\in\BZ$ one has
\[
\HH {-i}{\Hom RPQ}=\CH{i}{\Hom RPQ}\cong\CH{i}{\Hom RPN}=\Ext iRMN
   \]
in view of properties discussed in \ref{chunk:semiproj} and
\ref{chunk:cohomology}.  We set
   \[
\Ext *RMN=\bigoplus_{i\in\BZ}\Ext iRMN\,.
   \]
This is a graded module over $R^{\mathsf c}$, the center of the ring  
$R$.

Composition of homomorphisms turns $\Hom RQQ$ and $\Hom RPP$ into DG
algebras over the center $\cent R$ of $R$, and $\Hom RPQ$ into a left DG
module over the first and a right DG module over the second.  The  
actions
are compatible, so $\Ext *RNN$ and $\Ext *RMM$ become graded $R^{\mathsf
c}$-algebras and $\Ext *RMN$ a left-right graded bimodule over them.

These structures do not depend on choices of resolutions.
   \end{bfchunk}

\section{Cohomological supports}
\label{cohomological supports}

In this section $R$ denotes an associative ring.

\begin{bfchunk}{Cohomology operations.}
  \label{chunk:center}
A \textit{ring of central cohomology operations} is a commutative graded
ring $\mca$ equipped with a homomorphism of graded rings
\[
\zeta_M\col\mca\lra\Ext *RMM
\]
for each $M\in\dcat R$, such that for all $N\in\dcat R$ and $\xi\in 
\Ext *RMN$
one has
\begin{equation}
\label{eq:centrality}
\xi\cdot\zeta_M(a)=\zeta_N(a)\cdot\xi
\quad\text{for every}\quad a\in\mca\,.
\end{equation}
For $N=M$ this formula implies that $\zeta_M(\mca)$ is in the center of
$\Ext *RMM$.

We assume that $\mca$ is non-negatively graded and that $\mca^i=0$
for $i$ odd or $2\mca=0$; this hypothesis covers existing examples and
avoids sign trouble.
   \end{bfchunk}

\begin{bfchunk}{Scalars.}
  \label{chunk:scalars}
Using the standard identifications of rings
  \[
\Ext*RRR=\Hom RRR=R^{\mathsf o}\,,
\]
where $R^{\mathsf o}$ denotes the opposite ring of $R$, one sees from
\eqref{eq:centrality} that the homomorphism of rings $\zeta_R\col
\mca\to R^{\mathsf o}$ maps every element $a\in\mca^0$ to the center of
$R^{\mathsf o}$.  We identify the centers of $R^{\mathsf o}$ and $R$.
Formula \eqref{eq:centrality} then shows that the action of $a$ on
$\Ext*RMN$ coincides with the maps induced by left multiplication with
$\zeta_R(a)$ on $M$ or on $N$.
   \end{bfchunk}

For the next definition we use the notion of support introduced in
\ref{chunk:supp}.

\begin{bfchunk}{Cohomological supports.}
\label{chunk:hsup}
Let $\mca$ be a graded ring of central cohomology operations, as above.
For each pair $(M,N)$ of complexes we call the subset
   \[
\hhsup MN = \msup{(\Ext{*}RMN)}\subseteq \proj{\mca}
   \]
the \emph{cohomological support} of $(M,N)$.  The
cohomological support of $M$ is
   \[
\hsup M = \hhsup MM\,.
   \]
    \end{bfchunk}

The theorem below is the main result of this section.

\begin{theorem}
\label{thm:realizability}
Let $M$ and $N$ be complexes of $R$-modules.

If the graded $\mca$-module $\Ext *RMN$ is noetherian, then for every
closed subset $X$ of\, $\hhsup MN$ there exist complexes $M_X$ in $ 
\thick
RM$ and $N_X$ in $\thick RN$, such that the following equalities hold:
  \[
X=\hhsup{M_X}{N}= \hhsup{M_X}{N_X}=\hhsup{M}{N_X}\,.
  \]

Moreover, when $N=M$ one can take $N_X=M_X$.
  \end{theorem}

The proof appears at the end of the section.  Some of the preparatory
material is used repeatedly throughout the paper.

Let $d$ be an integer.  The $d$th \textit{shift} of a complex $M$
is the complex $\shift M$ with $(\shift^d M)_n=M_{n-d}$ for all $n$
and $\dd^{\shift^d M}=(-1)^d\dd^M$.  The $d$th \textit{twist} of
a graded $\mca$-module $\mcm$ is the graded module $\mcm(d)$ with
$\mcm(d)^j=\mcm^{d+j}$ for all $j$.

\begin{bfchunk}{Functoriality.}
   \label{chunk:functorial}
Let $M,M',M''$ and $N,N',N''$ be complexes of $R$-modules.

There exist canonical isomorphisms of graded $\mca$-modules:
\begin{gather}
  \label{eq:shift}
\Ext{*}R{\shift M}N(1)\cong\Ext{*}R{M}N\cong
\Ext{*}RM{\shift N}(-1)\,;
\\
  \label{eq:sum}
\begin{aligned}
  \Ext{*}R{M'\oplus M''}N &\cong\Ext{*}R{M'}N \oplus\Ext{*}R{M''}N\,;
  \\
\Ext{*}RM{N'\oplus N''} &\cong\Ext{*}RM{N'} \oplus\Ext{*}RM{N''}\,.
\end{aligned}
  \end{gather}
Indeed, basic properties of the functor $\Hom{\dcat R}--$ show that for
a fixed $N$ (respectively, $M$) the canonical isomorphisms of graded
$R^{\mathsf c}$-modules are linear for the action of $\Ext*RNN$ on the
left (respectively, of $\Ext*RMM$ on the right).  They are $\mca$-linear
because of the centrality of $\mca$, see \eqref{eq:centrality}.

Similarly, exact triangles $M'\to M\to M''\to$ and $N'\to N\to N''\to$
in $\dcat R$ induce exact sequences of graded $\mca$-modules
\begin{equation}
  \label{eq:sequence}
  \begin{gathered}
\xymatrixrowsep{.35pc}
\xymatrixcolsep{1pc}
\xymatrix{
\Ext{*}R{M''}N \ar@{->}[r] &\Ext{*}RMN \ar@{->}[r] &\Ext{*}R{M'}N \ar@ 
{->}[r] &
\\
\Ext{*}R{M''}N(1)\ar@{->}[r]^{\quad\quad} &\Ext{*}RMN(1)
  \\
\Ext{*}RM{N'} \ar@{->}[r]&\Ext{*}RM{N} \ar@{->}[r]&\Ext{*}RM{N''}\ar@ 
{->}[r] &
\\
\Ext{*}RM{N'}(1) \ar@{->}[r] &\Ext{*}RMN(1)
}
\end{gathered}
\end{equation}
\end{bfchunk}

Putting together the remarks in \ref{chunk:supp} and
\ref{chunk:functorial}, one gets:

\begin{lemma}
  \label{lem:hsup}
\pushQED{\qed}
In the notation of \emph{\ref{chunk:functorial}} the following  
statements hold.
\begin{gather}
  \label{eq:hsup-shift}
\hhsup {\shift M}N=\hhsup {M}N=\hhsup M{\shift N}\,.
  \\
  \label{eq:hsup-sum}
\begin{aligned}
  \hhsup {M'\oplus M''}N&=\hhsup {M'}N\, \cup\, \hhsup {M''}N\,.
  \\
\hhsup M{N'\oplus N''}&=\hhsup M{N'}\, \cup\, \hhsup M{N''}\,.
\end{aligned}
  \\
\label{eq:hsup-sequence}
\begin{aligned}
\hhsup MN&\subseteq\hhsup {M'}N\, \cup\, \hhsup {M''}N\,.
  \\
\hhsup MN&\subseteq\hhsup M{N'}\, \cup\, \hhsup M{N''}\,.
\end{aligned}
\end{gather}

If $\Ext{*}RMN$ is eventually zero, then $\hhsup MN =\varnothing$.
The converse holds when $\Ext{*}RMN$ is eventually
noetherian over $\mca$.
   \qed
  \end{lemma}

The exact sequences \eqref{eq:sequence} imply the following statement:

\begin{lemma}
   \label{lem:noetherian}
Let $M$ be a complex of $R$-modules.

The full subcategory of $\dcat R$ consisting of complexes $L$ with $\Ext
*RM{L}$ (respectively, $\Ext *R{L}M$) eventually noetherian over $\mca$
is thick.
  \qed
  \end{lemma}

\begin{bfchunk}{Mapping cone.}
   \label{chunk:cone}
Let $M$ be a complex of $R$-modules.

For each $\vf\in\mca^d$ the morphism $\zeta_M(\vf)\col M\to{\shift^dM}$
defines an exact triangle
  \begin{equation}
  \label{eq:cone}
M\xra{\ \zeta_M(\vf)\ }\shift^dM\lra \koz \vf M\lra
  \end{equation}
which is unique up to isomorphism.

Let $N$ be a complex of $R$-modules and set $\mcm = \Ext *RMN$.
By \eqref{eq:sequence} and \eqref{eq:shift}, the triangle above yields
an exact sequence of graded $\mca$-modules
\begin{equation}
  \label{eq:sequence5}
\mcm(-d-1)\lra\mcm(-1)\lra \Ext *R{\koz\vf M}N\lra\mcm(-d)\lra\mcm\,.
\end{equation}
The maps at both ends are given by multiplication with $\vf$, so from
\eqref{eq:sequence5} one can extract an exact sequence of graded $\mca 
$-modules
\begin{equation}
  \label{eq:sequence3}
0\lra\big(\mcm/{\mcm\vf}\big)(-1)
   \lra \Ext *R{\koz \vf M}N  \lra (0:_{\mcm} \vf )(-d)\lra 0\,.
\end{equation}

Let $\bs{\vf}=\vf_1,\dots,\vf_n$ be homogeneous elements
in $\mca$.  Set $\bs{\vf}'=\vf_1,\dots,\vf_{n-1}$.  A complex
$\koz{\vf_n}{(\koz{\bs{\vf}'}M)}$ is defined uniquely up to isomorphism
in $\dcat R$;  we let $\koz{\bs{\vf}}M$ denote any such complex.
Iterated references to \eqref{eq:cone} yield
   \begin{equation}
  \label{eq:thick}
\koz{\bs{\vf}}M\in\thick RM\,.
   \end{equation}
  \end{bfchunk}

\begin{example}
\label{ex:cone}
If $\vf_1,\dots,\vf_n$ are in $\mca^0$, then $\zeta_M({\vf_i})$
is the homothety $M\to M$ defined by the central element
$z_i=\zeta_R(\vf_i)\in R$; see \ref{chunk:scalars}.  Thus, in $\dcat R$
one has $\koz{\bs\vf}M\simeq M\otimes_{\cent R}K(\bsz)$, where $K(\bsz)$
is the Koszul complex on $\bsz=z_1,\dots,z_n$.
  \end{example}

\begin{proposition}
\label{prop:cone}
Let $M$, $N$ be complexes of $R$-modules and $\bs{\vf}=\vf_1,\dots, 
\vf_n$
a sequence of homogeneous elements in $\mca$.

If the $\mca$-module $\Ext *RMN$ is eventually noetherian, then so  
are the
$\mca$-modules $\Ext *R{\koz{\bs{\vf}}M}N$, $\Ext *RM{\koz{\bs{\vf}}N}$,
and $\Ext *R{\koz{\bs{\vf}}M}{\koz{\bs{\vf}}N}$, and
   \begin{align*}
\hhsup {\koz{\bs{\vf}} M}N
& = \hhsup {\koz \vf M}{\koz{\bs{\vf}} N} = \hhsup M{\koz{\bs{\vf}} N}
  \\
&=\hhsup MN\, \cap\,  \msup (\mca/\mca{\bs{\vf}}) \,.
\end{align*}
  \end{proposition}

\begin{proof}
It suffices to treat the case when $\bs{\vf}$ has a single element,
$\vf$.  {}From the exact sequence \eqref{eq:sequence3} one sees that
$\Ext *R{\koz{\vf}M}N$ is eventually noetherian.

Set $\mcm=\Ext *RMN$.  The inclusion below holds because $(0:_{\mcm}
\vf )$ is a submodule of  $\mcm$ and a module over $\mca/\mca\vf$;
the equality comes from \ref{chunk:supp}(5):
\begin{align*}
\msup (0:_{\mcm} \vf ) &\subseteq \msup{\mcm}\, \cap\, \msup (\mca/ 
\mca \vf)\,;
  \\
\msup\big({\mcm}/{\mcm\vf}\big)
   & = \msup{\mcm}\, \cap\,  \msup (\mca/\mca \vf)\,.
\end{align*}
The exact sequence \eqref{eq:sequence3} and \ref{chunk:supp}(2) now  
imply
an equality
  \[
\hhsup {\koz \vf M}N = \msup{\mcm}\, \cap\,  \msup (\mca/\mca \vf)\,.
  \]
By a similar argument, $\Ext *RM{\koz{\vf}N}$ is eventually noetherian
and one has
   \[
\hhsup M{\koz \vf N} = \msup{\mcm}\, \cap\,  \msup (\mca/\mca \vf)\,.
   \]
The remaining equality is a formal consequence of those already  
available.
  \end{proof}

   \begin{proof}[Proof of Theorem~\emph{\ref{thm:realizability}}]
Set $\mcm=\Ext *RMN$ and $\mci=\ann_{\mca}\mcm$.

> From \ref{chunk:supp}(3) and \ref{chunk:supp}(4) we get
  \[
\hhsup MN = \msup \mcm =\msup(\mca/\mci)
  \]
As $\mcm$ is noetherian the closed subset $X$ of $\msup(\mca/\mci)$
has the form
  \[
X=\msup(\mca/\mci)\cap\msup(\mca/\mca\bs{\vf})
  \]
where $\bs{\vf}$ is a finite set of homogeneous elements of $\mca$;
see \ref{chunk:noetherian}.  Thus, one gets
  \[
X=\hhsup MN\cap\msup(\mca/\mca\bs{\vf})\,.
  \]
Choose complexes $M_X$ and $N_X$ representing $\koz{\bs{\vf}}M$ and
$\koz{\bs{\vf}}N$, respectively.  One has $M_X\in\thick R{M}$ and
$N_X\in\thick R{N}$, see \eqref{eq:thick}.  Also, one gets
   \[
X=\hhsup{M_X}N=\hhsup{M_X}{N_X}=\hhsup M{N_X}
   \]
from Proposition \ref{prop:cone}.  Clearly, when $M=N$ one may choose
$N_X=M_X$.
   \end{proof}

\section{Realizability by modules}
\label{realizability by modules}

In this section $R$ is an associative ring and $\mca$ is a ring of  
central
cohomology operations on $\dcat R$, see \ref{chunk:center}. The  
principal
result here is a partial enhancement of Theorem~\ref{thm:realizability}.
It contains Theorem~\ref{itheorem1} from the introduction.

\begin{etheorem}
\label{thm:realizability-modules}
Let $R$ be a left noetherian ring.

When $M$ and $N$ are complexes of $R$-modules with $\hh M$ and $\hh N$
finite, and $X$ is a closed subset of $\hhsup MN$ the following hold.
\begin{enumerate}[\quad\rm(1)]
   \item
If\, $\Ext*RMN$ is eventually noetherian over $\mca$ (and $X$ is
irreducible), then there exists a finite (and indecomposable) $R$-module
$M_X$ with
   \[
X=\hhsup {M_X}N
\quad\text{and}\quad
M_X\in\thick R{M\oplus R}\,.
   \]
\item
If, furthermore, $\Ext*RMR$ is eventually zero (and $X$ is irreducible),
then there exists a finite (and indecomposable) $R$-module $N_X$ with
   \[
X=\hhsup M{N_X}
\quad\text{and}\quad
N_X\in\thick R{N\oplus R}\,.
   \]
When $N=M$ one may choose $N_X=M_X$.
   \end{enumerate}
    \end{etheorem}

\begin{proof}
Using Theorem~\ref{thm:realizability}, choose complexes $C$ in $\thick
RM$ and $D$ in $\thick RN$, satisfying $\hhsup CN=X=\hhsup MD$.

By Lemma~\ref{lem:thick}(1), the $R$-modules $\hh C$ and $\hh D$ are
noetherian, so one has $\HH{\ges n}C=0=\HH{\ges n}D$ for some $n$.
Lemma~\ref{lem:thick} provides syzygy complexes $\syz nC$ in $\thick
R{M\oplus R}$ and $\syz nD$ in $\thick R{N\oplus R}$.  Another  
application
of Lemma~\ref{lem:thick}(1) shows that the following $R$-modules are
finite:
   \[
{M_X}=\HH 0{\syz n C} \quad \text{and} \quad  {N_X}=\HH 0{\syz nD}\,.
   \]
Lemma~\ref{lem:syzygies} yields $\syz nC\simeq M_X$ and $\syz nD\simeq
N_X$.

(1) comes from the equalities below, the second one given by
Lemma~\ref{lem:syzygies-supp}(1):
   \[
X=\hhsup CN=\hhsup {\syz n C}N=\hhsup {M_X}N\,.
   \]

(2)  As $\Ext *RMR$ is eventually zero, so is $\Ext
*R{M}F$; see Lemma~\ref{lem:thick}(3).  Thus, referring to
Lemma~\ref{lem:syzygies-supp}(2) for the second equality, one obtains
   \[
X=\hhsup MD=\hhsup M{\syz n D} = \hhsup M{N_X}\,.
   \]

When $N=M$ one can choose $D=C$ by Theorem~\ref{thm:realizability},
and hence get
   \[
N_X=M_X\simeq\syz nC\in\thick R{M\oplus R}\,.
   \]
Lemma~\ref{lem:thick}(3) now shows that $\Ext*R{M_X}F$ is eventually
zero when $F$ is free.  Thus, the already established assertion of the
theorem apply to $M_X$ and give
   \[
X=\hhsup {M_X}{M_X}\,.
   \]

It remains to establish the additional property when $X$ is
irreducible.  Being a noetherian module, $M_X$ is a finite direct sum
of indecomposables.  It follows from \eqref{eq:hsup-sum} that one can
replace $M_X$ with such a summand, without changing $\hhsup {M_X}N$.
A similar argument works for $N_X$.
   \end{proof}

The following general property of syzygy complexes was used above.

\begin{lemma}
  \label{lem:syzygies-supp}
Let $M$, $N$ be complexes of $R$-modules with bounded homology.

For every integer  $n$ the following hold.
\begin{enumerate}[\quad\rm(1)]
   \item
There is an equality $\hhsup MN = \hhsup {\syz n M}N$.
   \item
If $\Ext *RMF$ is eventually zero for every free module $F$, then also
   \[\hhsup MN = \hhsup M{\syz n N}\,.
  \]
   \end{enumerate}
    \end{lemma}

  \begin{proof}
(2)  Replacing $N$ with a semiprojective resolution, we may assume that
each $N_j$ is projective and $N_j=0$ for all $j\ll0$.  An elementary
argument using the formulas in \ref{chunk:functorial} shows that the
complexes $C$ of $R$-modules, for which $\Ext *RMC$ is eventually zero,
form a thick subcategory of $\dcat R$.  It contains the free modules by
hypothesis, and hence it contains all bounded complexes of projective
modules. Therefore, $\Ext{\ges i}RM{N_{<n}}=0$ holds for all $i\gg0$.

The inclusion $N_{<n}\subseteq N$ gives rise to an exact triangle
$N_{<n} \to N \to N_{\ge n}\to$ in $\dcat R$.  Again by \eqref 
{eq:sequence},
it induces an exact sequence of graded $\mca$-modules
   \[
\xymatrixrowsep{.35pc}
\xymatrixcolsep{1pc}
\xymatrix{
\Ext *RM{N_{<n}} \ar@{->}[r] &\Ext *RMN \ar@{->}[r] & \Ext *RM{N_{\ge  
n}} \ar@{->}[r] & \\
\Ext *RM{N_{<n}}(1) \ar@{->}[r] &\Ext *RMN(1)
}
  \]
In view of the preceding discussion it yields $\Ext{\ges
i}RMN\cong\Ext{\ges i}RM{N_{\ges n}}$ for all $i\gg 0$. On the other  
hand,
one has $N_{\ges n}=\shift^{n}\syz nN$ because $N$ is semiprojective.
The desired equality now follows from \ref{chunk:supp}(2).

(1)  This follows from a similar, and simpler, argument.
   \end{proof}

\section{Bialgebras and Hopf algebras}
\label{bialgebras}

In this section $k$ denotes a field.  We recall some notions concerning
bialgebras and Hopf algebras, referring to \cite{Mg} for details.

A \emph{bialgebra} over $k$ is a $k$-algebra $R$ with structure map
$\eta\col k\to R$ and product $\mu\col R\otimes_kR\to k$, equipped with
homomorphisms of rings $\eps\col R\to k$, the \textit{augmentation}
and $\Delta\col R\to R\otimes_kR$, the \textit{co-product}, satisfying
equalities
\begin{gather*}
\eps\eta=\idmap^k\,,\qquad
(\Delta\otimes \idmap^R)\Delta = \Delta(\Delta\otimes \idmap^R)\,,\\
\mu(\idmap^R\otimes \eta\eps)\Delta
=\idmap^R=\mu(\eta\eps\otimes \idmap^R)\Delta \,.
\end{gather*}
Given $R$-modules $M$, $N$ over a bialgebra $R$, the natural
$R\otimes_kR$-module structure on $M\otimes_kN$ restricts along
$\Delta$ to produce a canonical $R$-module structure:
\[
r\cdot (m\otimes n) = \sum_{i=1}^n (r_i'm\otimes_kr_i''n)
\quad\text{when}\quad
\Delta(r) =\sum_{i=1}^n (r_i'\otimes_kr_i'')\,.
\]
This extends to tensor products of complexes of $R$-modules.  Let $M$
be such a complex. The canonical isomorphisms below are easily seen to
be $R$-linear:
\begin{equation}
\label{eq:hopf-tensors}
k\otimes_k M\cong M
   \quad\text{and}\quad
M\otimes_k k\cong M\,.
\end{equation}

\begin{bfchunk}{Cohomology operations.}
   \label{hopf-ext}
Let $R$ be a bialgebra over $k$, and view $k$ as an $R$-module via
the augmentation $\eps$.  The ring $\Ext *Rkk$ has $\Ext0Rkk=k$ and is
\emph{graded-commutative}: for all $\alpha\in\Ext iRkk$ and $\beta\in 
\Ext
jRkk$ one has
   \[
\alpha\cdot \beta =(-1)^{ij}\beta \cdot\alpha\,;
   \]
see \cite[(VIII.4.7), (VIII.4.3)]{Mc} or \cite[(5.5)]{Msri}.  Thus,
every graded subring
\[
\mca\subseteq
\Ext{\scriptscriptstyle\bullet}Rkk=
\begin{cases}
\bigoplus_{i\ges 0} \Ext{2i}Rkk&\text{if }\operatorname{char}(k)\ne2 
\,;\\
\bigoplus_{i\ges 0} \Ext{i}Rkk&\text{if }\operatorname{char}(k)=2\,.
\end{cases}
  \tag{a}
  \]
is commutative.  The functor $-\otimes_kM$ preserves quasi- 
isomorphisms of
complexes of $R$-modules, so it induces a functor $-\otimes_kM\col 
\dcat R\to
\dcat R$.  In view of the isomorphism $k\otimes_kM\cong M$, see
\eqref{eq:hopf-tensors}, for each $M$ one gets a map
\[
\zeta_M\col\Ext{\scriptscriptstyle\bullet}Rkk\to \Ext *RMM\,.
\tag{b}
\]
It is readily verified to be a central homomorphism of graded $k$- 
algebras.
   \end{bfchunk}

The results in Section~\ref{cohomological supports} apply to any algebra
$\mca$ as above.  More comprehensive information is available for  
special
classes of bialgebras.

A \emph{Hopf algebra} is a bialgebra $R$ with a $k$-linear map
$\sigma\col R\to R$, the \emph{antipode}, satisfying $\eps\sigma=\eps$
and $\mu(1\otimes \sigma)\Delta = \mu(\sigma\otimes 1)\Delta $.  Quantum
groups offer prime examples.  A Hopf algebra is \emph{cocommutative}
if $\tau\Delta=\Delta$ holds, where $\tau(r\otimes s)=s\otimes r$.
For instance, for a group $G$ the $k$-linear maps defined by
\[
\eps(g)=1\,,
\quad
\Delta(g) = g\otimes g\,,
\quad\text{and}\quad
\sigma(g) = g^{-1}\quad
\text{for}\quad g\in G
\]
turn the group algebra $kG$ into a cocommutative Hopf algebra.  Other
classical examples are universal enveloping algebras of Lie
algebras and restricted universal enveloping algebras of $p$-Lie
algebras, where $p=\operatorname{char}(k)>0$.

\begin{bfchunk}{Finiteness.}
\label{hopf:finiteness}
Let $R$ be a Hopf algebra such that $\rank_k R$ finite.

If $R$ is cocommutative, then $\Ext{\scriptscriptstyle\bullet}Rkk$
is finitely generated as a $k$-algebra, and $\Ext*RMN$ is a finite
$\Ext{\scriptscriptstyle\bullet}Rkk$-module for all $R$-modules
$M,N$ of finite $k$-rank: This is a celebrated theorem of Friedlander
and Suslin~\cite[(1.5.2)]{FS}, which extends earlier results for group
algebras (Evens, Golod, Venkov)  and for restricted Lie algebras
(Friedlander and Parshall).

It is not known whether cohomology has similar finiteness properties
when $R$ is not cocommutative; for positive solutions in interesting
classes of such Hopf algebras see Pevtsova and Witherspoon~\cite{PW},
and the bibliography there.
  \end{bfchunk}

\begin{bfchunk}{Cohomological varieties.}
   \label{hopf:variety}
Let $R$ be a Hopf algebra with $\rank_k R$ finite, set $\mca=
\Ext{\scriptscriptstyle\bullet}Rkk$, see \ref{hopf-ext}(a), and let
$\ov k$ be an algebraic closure of $k$.

For a complex $M$ with $\Ext *RMM$ eventually noetherian over $\mca$
define (with notation as in \ref{chunk:var}) the \emph{cohomological
variety} of $M$ to be the subset
   \[
\hhvar[R]M = \mvar {\Ext *RMM} \subseteq \Max{(\mca\otimes_k\ov k)} \,.
   \]
    \end{bfchunk}

\begin{etheorem}
\label{hopf:main}
Let $R$ be a Hopf algebra over $k$, such that $\rank_kR$ is finite;
set $\mca=\Ext{\scriptscriptstyle\bullet}Rkk$.  Let $M$ be a complex
with $\hh M$ finite over $R$.

If $\Ext*RMM$ is eventually noetherian over $\mca$, then for each closed
conical $k$-rational subset $X$ of $\hhvar M$ there is a finite $R$- 
module
$M_X$, such that
   \[
X=\hhvar {M_X}
   \quad\text{and}\quad
M_X\in\thick R{M\oplus R}\,.
   \] \end{etheorem}

\begin{proof}
Hopf algebras of finite rank are self-injective, see
\cite[(2.1.3)(4)]{Mg}, so one has $\Ext{\ges1}R-R=0$.  It remains
to invoke Theorem~\ref{thm:realizability-modules} and refer to
\ref{chunk:var}.
  \end{proof}

  \begin{bfchunk}{Applications.}
In view of \ref{hopf:finiteness}, for $M=k$ the theorem specializes
to results of Carlson~\cite{Ca}, Suslin, Friedlander, and
Bendel~\cite[(7.5)]{SFB}, Pevtsova and Witherspoon~\cite[(4.5)]{PW},
among others.  It is clear that there are also versions dealing with
supports of pairs of modules, and with complexes.
  \end{bfchunk}

\section{Associative algebras}
\label{Associative algebras}

Here $k$ is a field and $R$ is a $k$-algebra.  Let $R^{\mathsf{o}}$  
denote
the opposite algebra of $R$, set $R^\sse = R\otimes_k R^{\mathsf{o}}$,
and turn $R$ into a left $R^\sse$-module by $(r\otimes r')\cdot s=  
rsr'$.

\begin{bfchunk}{Cohomology operations.}
   \label{ass-operations}
The \emph{Hochschild cohomology} of $R$ is the $k$-algebra
  \[
\hch *kR =  \bigoplus_{i\ges 0}\Ext {i}{R^\sse}RR\,.
  \]
Gerstenhaber~\cite[Cor.~1]{Ge} proved that it is graded-commutative,
so any subring
  \[
\mca\subseteq\hch{\scriptscriptstyle\bullet}kR=\bigoplus_{i\ges
0}\hch{2i}kR
   \tag{a}
  \]
is commutative.  The map $r\mapsto 1\otimes r$ is a homomorphism of  
rings
$R^{\mathsf{o}}\to R^\sse$.  It turns each complex of $R^\sse$-modules
into one of right $R$-modules.  Thus, $-\otimes_RM$ is an additive
functor from complexes of $R^\sse$-modules to complexes of $R$-modules,
where $R$ acts on the target via the homomorphism of rings $R\to R^\sse$
given by $r\mapsto r\otimes 1$.  It induces an exact functor $-\lotimes
RM\col \dcat{R^\sse}\to \dcat R$ of derived categories, which produces
a homomorphism of graded rings
  \[
\Ext *{R^\sse}RR\to \Ext *R{R\lotimes RM}{R\lotimes RM}\,.
  \]
The isomorphism $R\lotimes RM\simeq M$ now yields a natural homomorphism
  \[
\zeta_M\col \hch{\scriptscriptstyle\bullet}kR \to \Ext *RMM \tag{b}
  \]
of graded rings. These maps satisfy condition \eqref{eq:centrality},
see \cite[(10.1)]{So}.
  \end{bfchunk}

The results in Sections \ref{cohomological supports} and
\ref{realizability by modules} apply to any algebra $\mca$ as above.  
Once
again, we focus on a special case to relate them to available  
literature.

   \begin{bfchunk}{Finiteness.}
   \label{ass:finiteness}
Let $R$ be a $k$-algebra with $\rank_kR$ finite, $J$ the Jacobson
radical of $R$, and set $K=R/J$.  It is rarely the case that the $\hch
*kR$-module $\Ext *RM{K}$ is noetherian for every finite $R$-module $M$;
see \cite[\S1]{EHSST}.  Examples when this property holds include the  
Hopf
algebras in \eqref{hopf:finiteness}, exterior algebras, and commutative
complete intersections rings, see \eqref{lci-operations}.
   \end{bfchunk}

\begin{bfchunk}{Cohomological varieties.}
\label{ass:support}
For $R$ as in \eqref{ass:finiteness}, let $\mca$ be a subring
of $\hch{\scriptscriptstyle\bullet}kR$ with $\mca^0=k$, see
\ref{ass-operations}, and let $\ov k$ be an algebraic closure of $k$.

For a complex of $R$-modules $M$, such that the $\mca$-module $\Ext
*RM{K}$ is eventually noetherian, define the \emph{cohomological  
variety}
of $M$ to be the subset
   \[
\hhvar M=\mvar {\Ext *RM{K}}\subseteq \Max{(\mca\otimes_k\ov k)}\,.
   \]
As $\mca$ acts on $\Ext *RM{K}$ through $\Ext *R{K}{K}$, one has $\hhvar
M\subseteq\hhvar K$.
   \end{bfchunk}

\begin{etheorem}
\label{ass:main}
Let $R$ and $\mca$ be as in \emph{\ref{ass:support}}, and let $M$ be a
complex of $R$-modules with $\hh M$ finite over $R$.

If\,\ $\Ext *RMK$ is eventually noetherian over $\mca$, then for each
closed conical $k$-rational subset $X$ of $\hhvar M$ there is a finite
$R$-module $M_X$, such that
   \[
X=\hhvar {M_X}
   \quad\text{and}\quad
M_X\in\thick R{M\oplus R}\,.
   \]
  \end{etheorem}

  \begin{proof}
The $R$-module $R$ admits a finite filtration with subquotients  
isomorphic
to direct summands of $K$.  Thus, when the $\mca$-module $\Ext *RM{K} 
$ is
noetherian, so is $\Ext *RMR$. On it $\mca$ acts through $\Ext *RRR= 
\cent
R$; see \ref{chunk:scalars}. This means $\Ext{\gg0}RMR=0$, so we may
use Theorem~\ref{thm:realizability-modules}, then \ref{chunk:var}.
  \end{proof}

  \begin{bfchunk}{Application.}
When the $\mca$-module $\Ext *R{K}{K}$ is noetherian,
Theorem~\ref{ass:main} with $M=K$ yields a result of Erdmann
\textit{et al}; see \cite[(3.4)]{EHSST}.
    \end{bfchunk}

\section{Commutative local rings} \label{local complete intersections}

We say that $(R,\fm,k)$ is a \textit{local ring} if $R$ is a commutative
noetherian ring with unique maximal ideal $\fm$ and $k=R/\fm$.

An \textit{embedded deformation} of codimension $c$ of $R$ is a  
surjective
homomorphism $\varkappa\col Q\to R$ of rings with $(Q,\fq,k)$ a local
ring and $\Ker(\varkappa)$ an ideal generated by a $Q$-regular sequence
in $\fq^2$, of length $c$.

\begin{bfchunk}{Cohomology operations.}
  \label{lci-operations}
Let $(R,\fm,k)$ be a local ring with an embedded deformation
$\varkappa\col Q\to R$ of codimension $c$.  Set
\[
\mca=R[\chi_1,\dots,\chi_c]
\tag{a}
   \]
where $\chi_1,\dots,\chi_c$ are indeterminates of degree $2$.  For each
$M\in\dcat R$ Avramov and Sun \cite[(2.7), p.~700]{AS} construct a
natural homomorphism of graded rings
\[
\zeta_M\col\mca\to\Ext{*}RMM \tag{b}
   \]
satisfying condition \eqref{eq:centrality}; when $M$ and $N$
are $R$-modules the resulting structure of graded $\mca$-module
on $\Ext*RMN$ coincides with that defined by Gulliksen \cite{Gu}.
For complexes $M$ and $N$, the $\mca$-module $\Ext{*}RMN$ is finite if
and only if the $Q$-module $\Ext{*}QMN$ is finite; see \cite[(5.1)]{AS}
and \cite[(4.2)]{AGP}.

The action of $\mca$ on $\Ext*RMk$ factors through the graded ring
\[
\mcr=\mca\otimes_Rk=k[\chi_1,\dots,\chi_c]\,.\tag{c}
\]
   \end{bfchunk}

Recall that a complex of $Q$-modules is said to be \emph{perfect} if  
it is
isomorphic, in $\dcat Q$, to a bounded complex of finite free $Q$- 
modules.

   \begin{lemma}
   \label{le:perfect}
Let $Q$, $R$, and $\mcr$ be the rings in \emph{\ref{lci-operations}}.

The following conditions are equivalent for each complex $M$ of $R$- 
modules:
   \begin{enumerate}[\quad\rm(i)]
\item $M$ is perfect over $Q$.
\item $\hh M$ is finite over $R$ and $\Ext *RMk$ is finite over $\mcr$.
   \end{enumerate}
    \end{lemma}

\begin{proof}
(i) $\implies$ (ii) As $M$ is perfect over $Q$, the $Q$-module $\Ext
*QMk$ is finite, and hence the $\mcr$-module $\Ext *RMk$ is finite;
see \ref{lci-operations}.

(ii) $\implies$ (i) It follows from \ref{lci-operations} that $\Ext
*QMk$ is finite over $Q$, and hence is eventually zero.  Since $\hh M$
is finite over $Q$, the complex $M$ admits a semiprojective  
resolution $F$
with each $F_i$ finite, $F_i=0$ for $i\ll 0$, and $\dd(F)\subseteq\fm  
F$;
see, for example, \cite{AFH}. This yields an isomorphism $\Ext iQMk\cong
\Hom Q{F_i}k$.  Thus $\Ext iQMk=0$ for $i\ge n$ implies $F_i=0$ for $i 
\ge
n$, so $F$ is a perfect complex of $Q$-modules that is quasi-isomorphic
to $M$.
  \end{proof}

\begin{bfchunk}{Cohomological varieties.}
   \label{lci-varieties}
Let $(R,\fm,k)$ be a local ring with an embedded deformation $\varkappa$
of codimension $c$, as in \ref{lci-operations}, and $\ov k$ an algebraic
closure of $k$.  Let $M$ be a complex of $R$-modules with $\Ext*{R}{M}k$
noetherian over $\mcr$. The \textit{cohomological variety} of $M$ is
the subset $\hhvar[\varkappa] M$ of ${\ov k}^{c}$ defined by the formula
   \[
\hhvar[\varkappa] M=\mvar[\mcr]{\Ext *RMk}\subseteq\Max{(\mcr\otimes_k 
\ov
k)} ={\ov k}^{c}\,,
   \]
where the second equality comes from Hilbert's Nullstellensatz.
When $M$ is a module and $\bsf$ is a $Q$-regular sequence that
generates $\Ker(\varkappa)$, the construction above yields the cone
$\hhvar[R]{\bsf;M}$ defined in \cite{Av:vpd}.
   \end{bfchunk}

A $Q$-module is a perfect complex in $\dcat Q$ if and only if $\pd_QM 
$ is
finite.  Thus, Theorem~\ref{itheorem2} from the introduction is obtained
from the next result by replacing affine cones by their  
projectivizations.

   \begin{etheorem}
   \label{thm:restated}
Let $(Q,\fq,k)$ be a local ring, ${\bsf}\subset \fq^2$ a $Q$-regular
sequence of length $c$, and $\varkappa\col Q\to Q/Q\bsf=R$ the canonical
surjection.

The assignment $ M\mapsto \hhvar[\varkappa]M$, which maps complexes
in $\dcat R$ that are perfect over $Q$ to closed $k$-rational cones in
${\ov k}^c$, has the following properties:
   \begin{enumerate}[\quad\rm(1)]
\item
   It is surjective, even when restricted to modules.
\item
   $\hhvar[\varkappa]M=\{0\}$ if and only if $M$ is perfect over $R$.
\item
$\hhvar[\varkappa] M=\hhvar[\varkappa]{\Omega^R_n(M)}$ for every  
syzygy complex
  $\Omega^R_n(M)$.
  \item
$\hhvar[\varkappa]M=\hhvar[\varkappa]{M/\bsx M}$ for a module  $M$
and an $M$-regular sequence $\bsx$.
    \end{enumerate}
  \end{etheorem}

We import some material for the proof of part (1) of the theorem.

  \begin{chunk}
  \label{lci-example}
Under the hypotheses of the theorem Avramov, Gasharov, and Peeva
\cite[(3.3), (3.11), (6.2)]{AGP} give a nonzero finite $R$-module $G$  
with
   \[
\Ext{*}RGk\cong\mcr\otimes_k\Ext{*}QGk\,,
  \]
which also satisfies the conditions $\pd_QG<\infty$ and
$\Ext{\ges1}RGR=0$.
   \end{chunk}

   \begin{proof}[Proof of Theorem~\emph{\ref{thm:restated}}]
(1) One has $\mvar[\mcr]{\mcr}={\ov k}^{c}$ by the Nullstellensatz. In
view of \ref{chunk:var} it thus suffices to  find a module $G$ with
$\hhsup[\mcr]Gk=\spec\mcr$, and to show that every closed subset $X$ of
$\hhsup[\mcr]Gk$ is realizable by a module $M_X$ with $\pd_QM_X<\infty$.

The $R$-module $G$ from \ref{lci-example} has the
necessary property, as $\Ext{*}RGk$ is a nonzero graded free
$\mcr$-module. Theorem~\ref{thm:realizability-modules} yields a module
$M_X$ with the desired cohomological support and is in $\thick {R}{G 
\oplus
R}$. Since $G\oplus R$ is perfect over $Q$, the last condition implies
that so does $M_X$.

(2) Evidently $\hhvar[\varkappa]M=\{0\}$ if and only if
$\Supp[\mcr]{\Ext*RMk}=\varnothing$.  As the $\mcr$-module $\Ext*RMk$ is
finite, this is equivalent to $\$\Ext{\gg0}RMk=0$. Lemma~\ref 
{le:perfect},
applied with $Q=R$, yields the desired equivalence.

(3) This follows from Lemma~\ref{lem:syzygies-supp}.

(4) As noted in Example~\ref{ex:cone}, the complex $\koz {\bsx}M$
is quasi-isomorphic to the Koszul complex on $\bsx$, and hence to
the $R$-module $M/\bsx M$. Thus, Proposition~\ref{prop:cone} implies
$\hhvar[\varkappa]{M/\bsx M} = \hhvar[\varkappa]M $; see \ref 
{chunk:var}.
   \end{proof}

In Theorem~\ref{thm:restated} the hypothesis that $R$ has an embedded
deformation can be weakened in a useful way.  The main property is (1),
so we focus on it.

\begin{bfchunk}{Completions.}
   \label{completions}
The $\fm$-adic completion of $(R,\fm,k)$ is a local
ring, $(\wh R,\wh\fm,k)$.  The maps $R\to\wh R$ and
$M\to\wh R\otimes_R M=\wh M$ induce isomorphisms
   \[
\Ext*{\wh R}{k}k\lra\Ext*{R}{k}k
  \quad\text{and}\quad
\Ext*{\wh R}{\wh M}k\lra\Ext*{R}{M}k\,,\tag{d}
   \]
the first one of graded $k$-algebras, the second of graded modules,
equivariant over the first.  Thus, when $\wh R$ has an embedded
deformation $\varkappa$ of codimension $c$ the ring $\mcr$ from
\ref{lci-operations} acts on $\Ext*{R}{M}k$ for each $M\in\dcat R$.
As in \ref{lci-varieties}, when $\Ext*{R}{M}k$ is noetherian over $\mcr$
we define a \textit{cohomological variety} by:
   \[
\hhvar[\varkappa] M=\mvar[\mcr]{\Ext *RMk}\subseteq{\ov k}^{c}\,.
   \]

Observe that if $\varkappa$ is an embedded deformation of $R$, then $\wh
\varkappa$, completion with respect to the maximal ideal of $Q$, is an
embedded deformation of $\wh R$, so (d) above yields $\hhvar[\varkappa]
M=\hhvar[\wh \varkappa]{\wh M}$.
   \end{bfchunk}

The following descent result is of independent interest.

\begin{theorem}
   \label{thm:lci-completion}
Let $(R,\fm,k)$ be a local ring, $Q\to R$ an embedded deformation,
and $L$ a complex of ${\wh R}$-modules.

If\,\ $L$ is perfect over $Q$, then there exists a finite $R$-module
$M$ with
   \[
\hhvar[\varkappa]{M}=\hhvar[\varkappa]{L}
   \quad\text{and}\quad
\pd_Q\wh{M}<\infty\,.
   \]
  \end{theorem}

   \begin{proof}
By \ref{chunk:var}, it suffices to prove
$\msup {\Ext *RMk}=\msup {\Ext *RLk}$.

Choose a set of generators of the ideal $\fm$ and let $\bsx$ be its  
image
under the composition $R\to\wh R\to\mca$.  The equality $\bsx\Ext *{\wh
R}{L}k=0$ implies an inclusion
\[
\msup{\Ext *{\wh R}Lk}\subseteq\msup{(\mca/(\bsx)\mca)}\,.
\]
This yields the second equality below; the first one holds by
Lemma~\ref{lem:noetherian}:
   \begin{align*} \msup{\Ext *{\wh R}{\koz {\bsx}L}k}
     & = \msup {\Ext *{\wh R}Lk} \cap
     \msup{(\mca/(\bsx)\mca)}\\
      & = \msup {\Ext *{\wh R}Lk}\,.
  \end{align*}

The complex $\koz {\bsx}L$ is quasi-isomorphic to the Koszul
complex on $\bsx$ with coefficients in $L$; see \ref{ex:cone}.  Thus,
$\hh{\koz{\bsx}L}$ has finite length over $\wh R$, hence also over $R$.
Let $F\to \koz{\bsx}L$ be a semi-projective resolution over $R$ with $F$
a finite free complex.  One then has quasi-somorphisms
   \[
\wh R\otimes_R F \simeq \wh R\otimes_R (\koz {\bsx}L) \simeq \koz  
{\bsx}L
  \]
due to the flatness of $\wh R$ over $R$ and, for the second one, also to
the finiteness of the length of $\hh{\koz {\bsx}L}$ over $R$.  Fix
$n$ so that $\HH{\ges n}{\koz {\bsx}L}=0$ holds and set $M=\HH n{F_{\ges
n}}$.  One then has quasi-isomorphisms of complexes of $\wh R$-modules
\[
\wh M\cong \wh R\otimes_R M \simeq \wh R\otimes_R(F_{\ges n})\,.
  \]
They imply $\wh M\cong\syz [{\wh R}]i{\koz {\bsx}L}$, hence the
first equality below:
  \begin{align*}
\msup {\Ext *{\wh R}{\wh M}k}
  & = \msup {\Ext *{\wh R}{\syz i{\koz{\bsx}L}}k}
\\
  & = \msup {\Ext *{\wh R}{\koz {\bsx}L}k}\,.
  \end{align*}
Lemma~\ref{lem:syzygies} gives the second one. It remains to note that
since $L$ and $\wh R$ are both perfect over $Q$, so is any complex in
$\thick {\wh R}{L\oplus \wh R}$. Thus, $\koz{\bsx}L$ is perfect over
$Q$, by \eqref{eq:thick}, and hence so is $\wh M= \syz [{\wh R}]i{\koz
{\bsx}L}$, by Lemma~\ref{lem:syzygies}.
   \end{proof}

\begin{etheorem}
   \label{thm:lci-variety}
Let $(R,\fm,k)$ be a local ring.

If $\varkappa\col Q\to \wh R$ is an embedded deformation of codimension
$c$, then for each closed $k$-rational cone $X\subseteq {\ov k}^{c}$
there exists a finite $R$-module $M$ with
   \[
X=\hhvar[\varkappa] {M}
   \quad\text{and}\quad
\pd_Q\wh{M}<\infty\,.
   \]
  \end{etheorem}

\begin{proof}
Theorem~\ref{thm:restated}(1) provides a finite $\wh R$-module $L$ with
$X=\hhvar[\varkappa] L$, and of finite projective dimension over $Q$.
Now apply Theorem~\ref{thm:lci-completion}.
    \end{proof}

   \begin{bfchunk}{Application.}
\label{chunk:lci}
A local ring $R$ is \emph{complete intersection} if $\wh R$ admits an
embedded deformation $\varkappa\col Q\to R$ where $Q$ is a regular local
ring; see \cite[\S 29]{Ma}.  For such an $R$ Theorem \ref{thm:lci-variety}
specializes to a result proved by Bergh~\cite[(2.3)]{Be}, who uses
Tate cohomology, and by Avramov and Jorgensen \cite{AJ}, who establish
an existence theorem for cohomology modules by using equivalences of
triangulated categories and Koszul duality.
   \end{bfchunk}

     \end{document}